# The Geometry of Trifocal Curves with Applications in Architecture, Urban and Spatial Planning

*Maja Petrović*, University of Belgrade, Faculty of Transport and Traffic Engineering, Belgrade, Serbia;
Address: Vojvode Stepe 305, 11000 Beograd, Serbia,
telephone: +381 11 3091 259, fax: +381 11 3096 704, e-mail: majapet@sf.bg.ac.rs
*Bojan Banjac*, University of Belgrade, Faculty of Electrical Engineering, Belgrade, Serbia,
University of Novi Sad, Faculty of technical sciences – Computer Graphics Chair, Novi Sad, Serbia
*Branko Malešević*, University of Belgrade, Faculty of Electrical Engineering, Belgrade, Serbia

In this paper we consider historical genesis of trifocal curve as an optimal curve for solving the Fermat's problem (minimizing the sum of distance of one point to three given points in the plane). Trifocal curves are basic plane geometric forms which appear in location problems. We also analyze algebraic equation of these curves and some of their applications in architecture, urbanism and spatial planning. The area and perimeter of trifocal curves are calculated using a Java application. The Java applet is developed for determining numerical value for the Fermat-Torricelli-Weber point and optimal curve with three foci, when starting points are given on an urban map. We also present an application of trifocal curves through the analysis of one specific solution in South Stream gas pipeline project.

## 1. Historical concerns of optimal location

The Fermat problem is given in original Latin as (Fermat, 1679): "*datis tribus punctis, quartum reperire, a quo si ducantur tres rectae ad data puncta, summa trium harum rectarum sit minima quantitas*" or in the English translation "for three given points, the fourth is to be found, from which if three straight lines are drawn to the given points, the sum of the three lengths is minimum" (Brazil *et al.*, 2013). French mathematician Pierre de Fermat in XVII century in a private letter to Italian physicist and mathematician Evangelista Torricelli (Torricelli *et al.*, 1919) stated a problem of finding a point with following property: the sum of distances from one point to the vertices of a given triangle is minimal. The first solution of this problem was given by Torricelli with three equilateral triangles over the sides of the initial triangle and three circumcircles of equilateral triangles. Then the point of intersection of these circumcircles is solution of Fermat problem and it's called Fermat-Torricelli point F (see Fig. 1). This point in literature is known as the fifth significant point of the triangle (Mladenović, 2004) and this point is said to be the first triangle center discovered after ancient Greek times (Volek, 2006).

Next step in the research of Fermat's problem was given by Bonaventura Cavalieri. He proved that three line segments from vertex of the triangle through the Fermat-Torricelli point (short dash line see Fig. 1) determined angles of 120° (Kirszenblat, 2011).

Vincenzo Viviani, the pupil of Torricelli and disciple of Galileo, had published Torricelli's solution (Pérgamo and Viviano, 1695), which in modern terminology was given by (Heinrich, 1965):
1. If all angles in $\Delta ABC$ are less than 120°, the point $F$ that minimizes the sum of its distances from vertices $A$, $B$, $C$ is the that point inside $\Delta ABC$ at which $\sphericalangle AFB = \sphericalangle BFC = \sphericalangle CFA = 120°$.
2. If some angle of $\Delta ABC$ is 120° or more, say the angle at $C$, then $F = C$.

One century later, English mathematician Thomas Simpson in his paper *Doctrine and Application of Fluxions* (Simpson, 1750) simplified the previous Torricelli construction. He connected the outside vertices $A_1$, $B_1$, and $C_1$ of equilateral triangles with appropriate vertices of the initial triangle $\Delta ABC$ (long dashed line see Fig. 1.). The point of intersection of Simpson's lines is the Fermat-Torricelli point F (Kirszenblat, 2011).

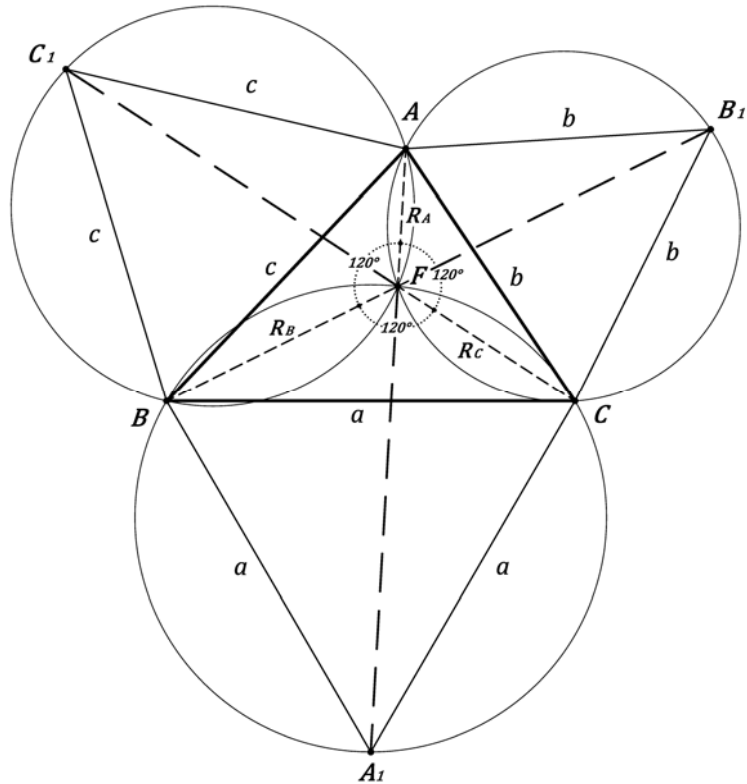

*Figure 1.The Torricelli construction of Fermat's point*

One upper bound of the sum of distance from the interior points to the vertices of initial triangle is given with the following inequality:

$$R_A + R_B + R_C < max\{a + b,\ b + c,\ c + a\}.$$

This inequality was given by Dutch mathematician Visschers (Crux, 1997) and it is proved by Romanian mathematician József Sándor (Sándor, 2005).

The Fermat-Torricelli point is the default point for optimal location in sense of the minimal sum equal weight distance from the point to vertices of triangle as considered by Austrian economist Alfred Weber (Weber, 1929). Weber's problem was later expanded on weighted problem with three points: find the location for factory such that transport costs to the factory from three resource suppliers are minimized, provided that the appropriate distances are multiplied by weight factors of objects that signify their importance to production process.

## 2. Location problems

These fundamentals of optimization problems of Fermat-Torricelli-Weber type are a good foundation for location analysis. Let us emphasize that Location analysis, Graph theory, Mathematical programming, Game theory and others belong to the area of Operational researches. Location problems themselves have great application in Urbanism, Architecture, Traffic, Industry (Boltyanski *et al.*, 1999, Drezner and Hamacher, 2002, Volek, 2006, Teodorović, 2009, Watanabe *et al.*, 2009). Current division of location problems to continuous, discrete and network (Mladenović, 2004) was based on the idea where to locate objects; if they are in the plane; what number of objects are serviced; in what way to allocate clients who demand service in objects … If it is possible to locate object in any point of considered region than the problem in question is continuous location problem. With discrete location problem object location is possible only in the predetermined set of points. With network model of the location problem, field of possible new location is anywhere on the given network; set of vertices on it is finite, and the set of arcs that connect vertices is continuum.

Classification of location problems (Teodorović, 2009) is also possible based on characteristics of certain location problems toward:

- Number of objects on network (one or more objects);
- Allowed places for object location (continuous or discrete location problems);
- Type of objects on network (medians, centers or objects with predefined system performance);
- Type of algorithm used for solving location problem (exact or heuristic algorithm);
- Number of criteria functions on which object locations are based (one or more criteria function).

### 3. Metric in location problems

Choice of a metric for measuring distance is of fundamental importance for solving location problems. Metric is defined by nature of the problem. Most common metric for solving location problems is Euclidean distance, and often that is the Euclidean distance with correction factor. Statistically it has been shown that measuring distances by roads gives 10-30% longer distances than corresponding Euclidean distance (Teodorović, 2009). If location problem in question is of urban nature and streets of the considered city are perpendicular to each other, than the distance is best determined by taxicab metric (Farahani and Hekmatflar, 2009). Location problems are most commonly connected to the plane. Still, as Earth's surface can be approximated by plane only in small dimensions, it is only logical that instead of Euclidean geometry, Riemann geometry is used. In that case the shortest distance between two cities (objects) is given by geodesic line. These lines are curves whose geodesic curvature is equal to zero in any point. On the sphere, geodesic line is great circle or orthodrome. Loxodrome is another geodesic line which is used for determining the distance between two points on sphere. Sea lines of communication are mostly loxodromic, so solving location problems in waterway traffic is conditioned also by loxodromic distance.

Let us note that Minkowski distance of the order $p$:

$$d_p = [|x_A - x_B|^p + |y_A - y_B|^p]^{\frac{1}{p}}$$

for two observed points $A(x_A, y_A)$ and $B(x_B, y_B)$ is considered for value of parameter $p \geq 1$ (Farahani and Hekmatflar, 2009). Especially for value of $p = 1$ Minkowski distance determines Manhattan, or the city block distance (short/long dashed lines see Fig. 2), while for $p = 2$ Minkowski distance determines Euclidian distance (continuous line see Fig. 2).

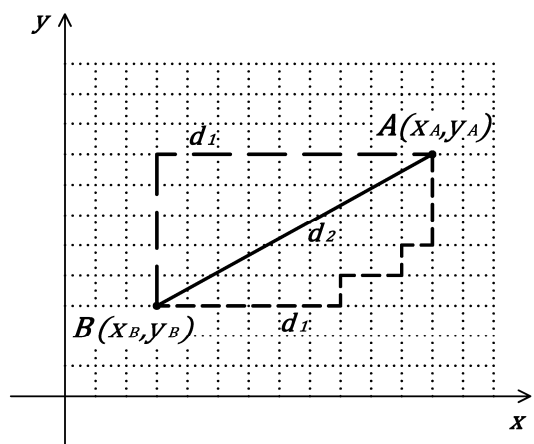

*Figure 2. Distances between two points*

Naturally, the term "distance" can also consider time distance, cost of travel or any other variable relative for location problem.

## 4. Trifocal curve as model for solving location problems

Scottish physicist and mathematician James Clerk Maxwell wrote his first scientific paper (Maxwell, 1847) in area of mechanical meaning of drawing mathematical curves (Harman, 1990): ellipse – curve with two focuses and curves with more than two focuses.

Standard definition of trifocal ellipse is given by the following equation:

$$R_A + R_B + R_C = S \qquad (1)$$

where are: $R_A = \sqrt{(x - x_A)^2 + (y - y_A)^2}$, $R_B = \sqrt{(x - x_B)^2 + (y - y_B)^2}$, $R_C = \sqrt{(x - x_C)^2 + (y - y_C)^2}$ Euclidean distances of the point $M(x, y)$ to three foci $A(x_A, y_A)$, $B(x_B, y_B)$, $C(x_C, y_C)$ and for given parameter $S > 0$. The Fermat-Torricelli point determines the minimal value of parameter $S = S_0$. If $0 < S < S_0$ then trifocal ellipse doesn't exist, and for $S > S_0$ trifocal ellipse is a non-degenerated egg curve. In that case the following names are applied for curve: egglipse, trisoid, polyellipse, trifocal ellipse, 3-ellipse, multifocal ellipse, etc (Sahadevan, 1974, 1987, Melzak and Forsyth, 1977, Erdös and Vincze, 1982, Sekino, 1999, Khilji, 2004).

These curves have oval shape and their construction is well known even from seventeenth century, but only in relation to Descartes construction of multifocal curves (Descartes, 1638). Maxwell construction has given certain simplifications (Mahon, 2003).

Weber's expansion of Fermat-Torricelli point takes in consideration the following equation

$$w_A R_A + w_B R_B + w_C R_C = S$$

where $w_A$, $w_B$ and $w_C$ are positive real values, weight factors.

Let us notice that researches of egg curves in engineering and their constructive geometrical solving are very contemporary even today because of ergonomics of their ovoid, bionic form (Rosin, 2000, 2004, Petrović, 2010, Barrallo, 2011). Also, the application of these curves on different problems of Fermat-Torricelli-Weber type are given in papers from theory of optimization (Soothill, 2010, Kupitz *et al.*, 2013).

P.V. Sahadevan was the first to introduce the term egglipse – a new curve with three focal points (Sahadevan, 1987). Sahadevan also introduced a classification of trifocal ellipses with collinear focuses in two categories. Trifocal ellipse is considered abnormal if focus $F_3$ is located on the middle of segment formed by $F_1$ and $F_2$, and if not, it is considered as normal (Sahadevan, 1974). By analyzing normal trifocal ellipse Sahadevan had found parametric equation of this curve and calculated area, using elliptic integral of first order. A. Varga and Cs. Vincze had also given similar parameterization of trifocal ellipse (1), (Varga and Vincze, 2008). The area of this closed curve was also determined using program package Maple. In literature we couldn't find procedure for finding perimeter and area of trifocal curve for non-collinear foci (Fig. 3) so special subroutine was developed in Java applet for calculation of their approximate numerical value.

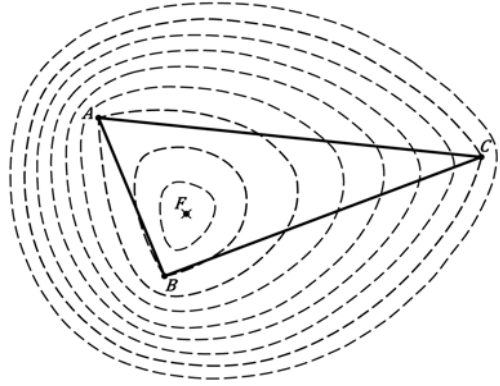

*Figure 3. Trifocal ellipses with non-collinear foci*

Subject of this paper is the trifocal ellipse (1). Let there be given expressions: $Q_1 = (x - x_A)^2 + (y - y_A)^2$, $Q_2 = (x - x_B)^2 + (y - y_B)^2$, $Q_3 = (x - x_C)^2 + (y - y_C)^2$ and $S = const, S > S_0$.

Let us consider the following algebraic transformations of equation (1):

$$\sqrt{Q_1} + \sqrt{Q_2} + \sqrt{Q_3} = S,$$

$$\sqrt{Q_1} + \sqrt{Q_2} = S - \sqrt{Q_3},$$

by first squaring we obtain

$$\left(\sqrt{Q_1} + \sqrt{Q_2}\right)^2 = \left(S - \sqrt{Q_3}\right)^2,$$

$$Q_1 + 2\sqrt{Q_1}\sqrt{Q_2} + Q_2 = S^2 - 2S\sqrt{Q_3} + Q_3,$$

$$2\sqrt{Q_1}\sqrt{Q_2} + 2S\sqrt{Q_3} = S^2 + Q_3 - Q_2 - Q_1,$$

by second squaring we obtain

$$\left(2\sqrt{Q_1}\sqrt{Q_2} + 2S\sqrt{Q_3}\right)^2 = (S^2 + Q_3 - Q_2 - Q_1)^2,$$

$$8S\sqrt{Q_1}\sqrt{Q_2}\sqrt{Q_3} = S^4 - 2S^2 Q_3 - 2S^2 Q_2 - 2S^2 Q_1 - 2Q_3 Q_2 - 2Q_3 Q_1 - 2Q_1 Q_2 + Q_1^2 + Q_2^2 + Q_3^2,$$

and by third squaring we finally obtain

$$\left(8S\sqrt{Q_1}\sqrt{Q_2}\sqrt{Q_3}\right)^2 = \left(S^4 - 2S^2 Q_3 - 2S^2 Q_2 - 2S^2 Q_1 - 2Q_3 Q_2 - 2Q_3 Q_1 - 2Q_1 Q_2 + Q_1^2 + Q_2^2 + Q_3^2\right)^2$$

which in the end determines algebraic curve of eight order:

$$64 S^2 Q_1 Q_2 Q_3 = ((S^2 + Q_3 - Q_2 - Q_1)^2 - 4Q_1 Q_2 - 4S^2 Q_3)^2.$$

The trifocal ellipse is one part of the previous curve, which is presented in (Nie *et al.*, 2008) using representation by the appropriate determinant. Analog conclusion also derives when considering Weber weighted variant of Fermat-Torricelli point. Let us emphasize that trifocal curves have natural application in multiple optimization problems. In the interior of trifocal curve of level $S$ we obtained the following inequality:

$$S \geq w_A R_A + w_B R_B + w_C R_C \geq S_0.$$

In the case $S = S_0$ we obtained Fermat-Torricelli-Weber point. If this point has inconvenient position in the sense of urban and spatial planning, then we determine the whole interior of trifocal curve as an appropriate solution for the point for which the sum of weighted distances is under or equal to the value of $S$.

## 5. Java application for visual representation of trifocal curves

In this paper we describe the application which was developed in programming language Java, whose purpose is to determine weighted Fermat-Torricelli point based on image input in application. Our application was developed in programming language Java since it was intended to be available over the internet. One of the basic goals of this project was for the application to be available for all researchers in this field, which asked for platform independent programming language, such as Java, to circumvent a need for several separate distributions of the application. Because the application was developed as an applet, and has a wide spectrum of functionalities, as well as new changes that company Oracle has introduced to runtime environment of programming language Java about digitally signing all files, when application is run, security question pops up which has to be answered confirmatively. As the application was intended to be available free of charge, no digital key from certified organization was available. Internet application is available as 3-Ellipses software on webpage http://symbolicalgebra.etf.bg.ac.rs/Java-Applications/3-Ellipses/index.html .

While developing Java applet, two types of graphical representations were implemented. The first type of graphical representation defines the interior of trifocal ellipse for the given triangle of foci and parameter $S$, as well as area and perimeter of region defined by it. Calculation of the area and the perimeter was determined by approximating graphical-numerical methods. Also accuracy estimate of considered calculation was determined. For higher precision it is necessary to use Calculus by computer algebra systems as Maple (Zuber and Štulić, 2006) or Mathematica (Groß and Strempel, 1998, Petruševski *et al.*, 2010). The second type of graphical representation defines trifocal isolines under condition $d_1 + d_2 + d_3 = S$ ($S \geq S_0$), which are represented in a graphic box. Parameter $S = S_0$ is value for which trifocal curve degenerates into Fermat-Torricelli point. Isolines are calculated only in defined rectangular region (graphic box), in which we also calculate points where the value of sum $d_1 + d_2 + d_3$ is minimal and maximal.

For any initial values of coordinates of the points *A*, *B* and *C*, the lengths $d_1$, $d_2$ and $d_3$ are modified into weighted lengths $w_A \cdot d_1$, $w_B \cdot d_2$, $w_C \cdot d_3$ by use of sliders where $w_A$, $w_B$, $w_C$ which have the range of decimal values from 0 to 10. Sliders $w_A, w_B, w_C$ have default value of 1 and are positioned in the right corner of the application main window. If graphical representation of curve is active, one additional slider is present for parameter $S$, which can have a decimal value from 0 to 2000.

In the applications menu bar several sections of program options are located. In section labeled Mode, users can select option Curve for graphical representation of trifocal ellipse or by selecting option Color mapping transfer to isolines representation (Fig. 4.a).

In menu labeled Opacity users can select a level of the opacity of graphical representations. Offered options are 20%, 40%, 60%, 80%, 100% and the custom level of the opacity where users are prompted in a dialog to input decimal value ranging from 0 to 100 (Fig. 4.b).

One of functionalities that was developed in the application is setting of the background map to graphical representation. In the menu labeled Background, users can select from several predefined maps for background, remove background image, or input background image of their choosing (Fig. 4.c).

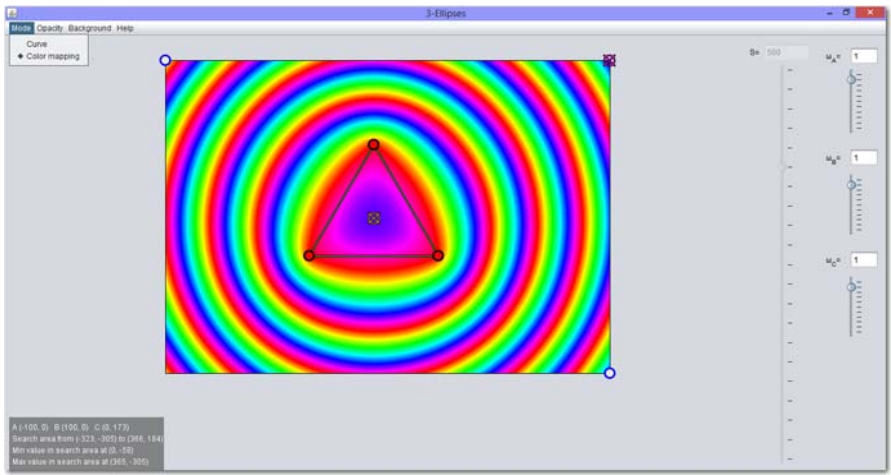

*a) Mode users*

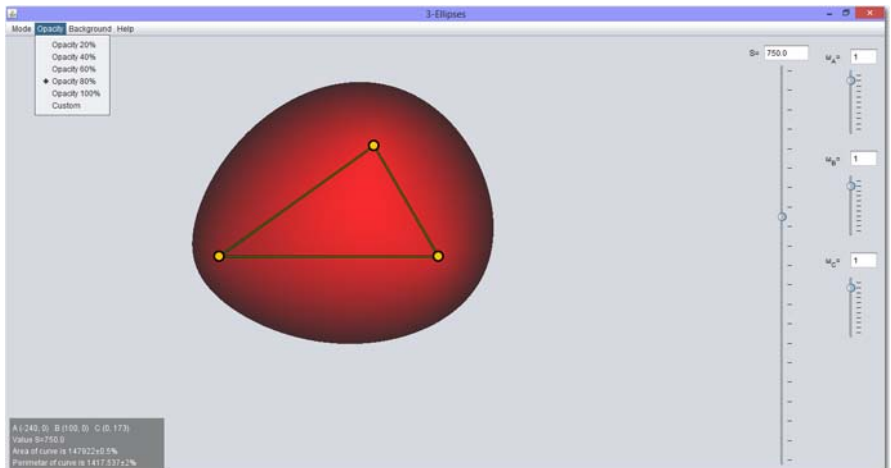

*b) Opacity*

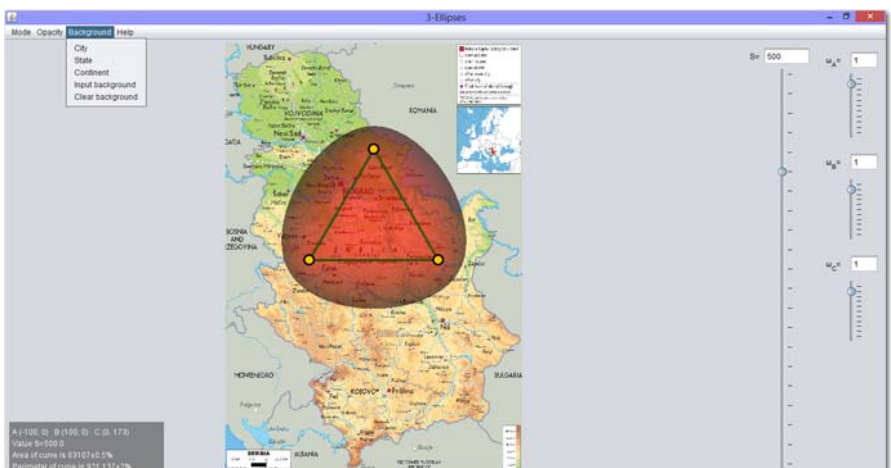

*c) Background*
*Figure 4. Java applet labels*

## 6. Examples of Application of Java Applet in Urban and Spatial Planning

In this part of the paper we shall conduct the analysis of application of Fermat-Torricelli point for three given points, i.e. locations/cities in whose near proximity South Stream pipeline is located.

Gas lines start from southern part of Russia, close to Beregovaya city at the bottom of the Black Sea, to Varna city in Bulgaria. The total length of the Black Sea section will exceed 900 kilometers. The onshore section in Europe will be 1,455 kilometers long (Fig. 5).

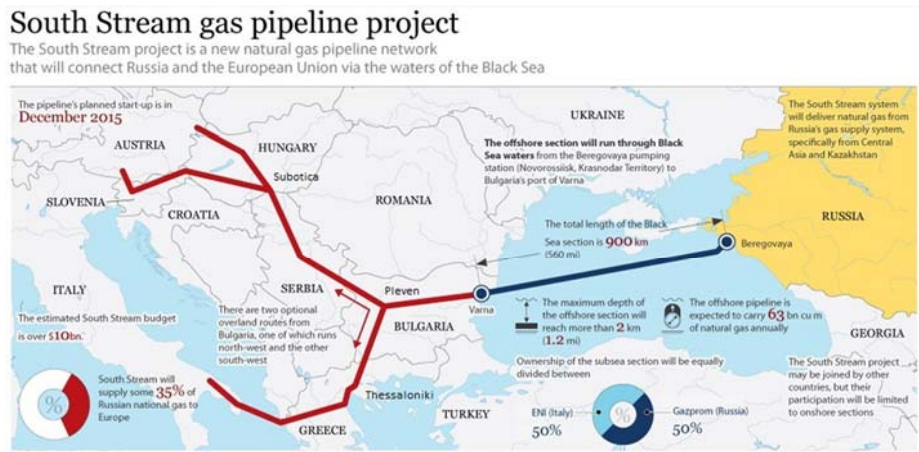

*Figure 5. South Stream gas pipeline project (Oikonomia, 2013)*

Near Pleven city in Bulgaria was planned that the pipeline would split in two ways. One part would run from Bulgaria, over the Greece, and under the water run to Italy, and other would run over Serbia and Hungary to Austria. South Stream will enter Serbia from Bulgaria near Zaječar city, eastern Serbia, and exit near Subotica city, northern Serbia.

If three given locations of Beregaya (Russia), Thessaloniki (Greece) and Subotica (Serbia) are observed, Fermat-Torricelli point would be near Pleven in Bulgaria (Fig. 6.*a* – see small trifocal curve around Pleven). On Figure 6.*b* near Fermat-Torricelli point isolines which are trifocal ellipses can be seen. Selecting Varna (Bulgaria) instead of the location Beregaya (Russia), the developed application determines that Fermat-Torricelli point also stays near city Pleven (Fig. 6. *c*). Let us notice that the optimal point in proximity of city Pleven is acquired with values of weight factors of one in this application. It infers that both sides of the pipelines to Thessaloniki and Subotica are equally important. Graphical representation of trifocal curve, calculation of its area and perimeter for set value of $S > S_0$ is available in the application. Let us note that similar web-based support of spatial planning in Serbia is considered in paper (Bazik and Dželebdžić, 2012)

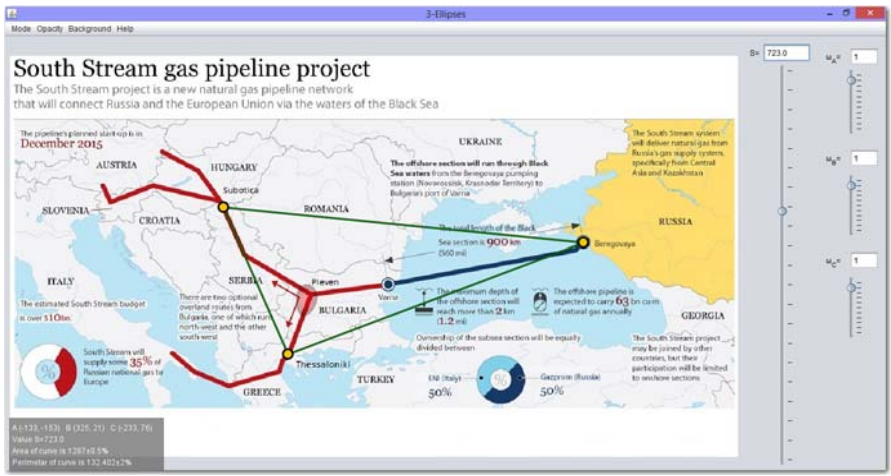

*a)*

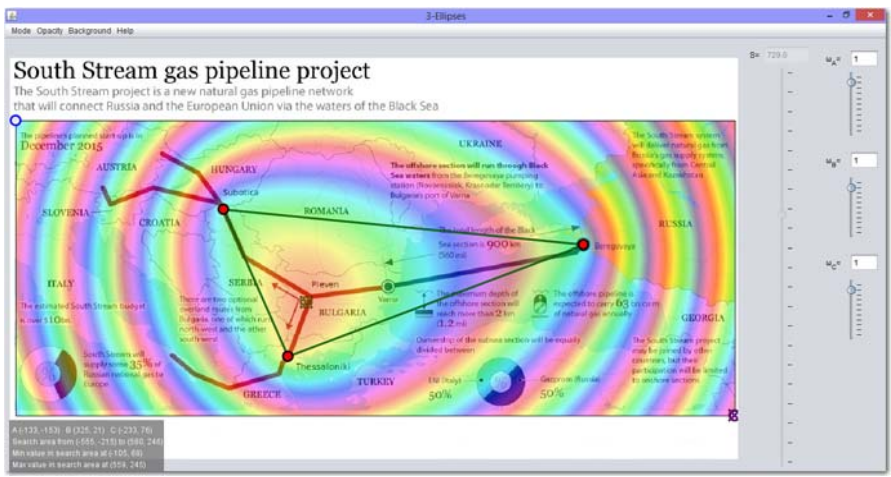

*b)*

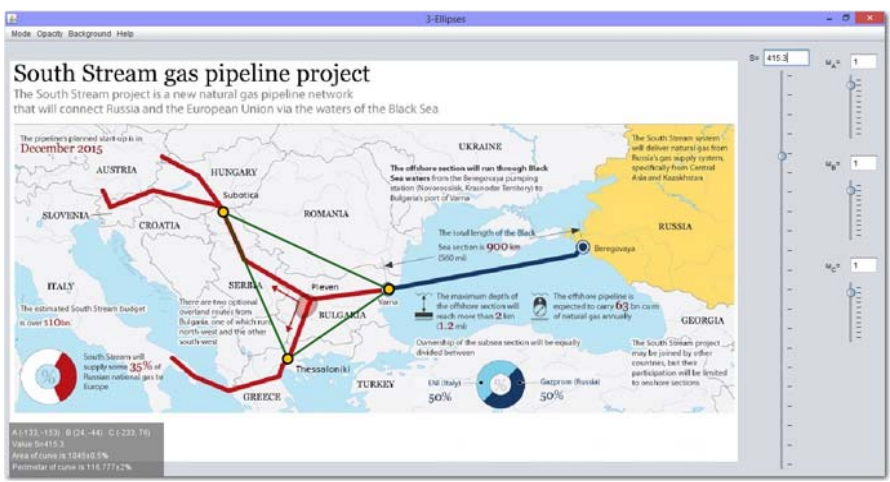

*c)*

*Figure 6. Application of Java applet*

## 7. Conclusion

The application of trifocal curves through analysis of one specific solution in urban and spatial planning was presented in this paper. The developed Java applet for plotting of these curves had enabled visual representation of optimal solution as well as fast calculation of numerical value of parameter $S$ (that is minimal sum of distances). The area and perimeter of trifocal curve for any initial values of coordinates of points $A$, $B$, $C$, parameter $S$, and weight factors $w_A$, $w_B$, $w_C$ were also developed.

Our research of optimal curves will further consider curves based on $n$ input points (where $n \geq 4$). The analytic and algebraic representation of these $n$-focal curves will also be considered. A similar line of research is also visual representation of these curves for other types of metric (Manhattan, Riemann or corrected Euclidean distance).

**Acknowledgement.** Research is partially supported by the Ministry of Science and Education of the Republic of Serbia, Grant No. III 44006 and ON 174032.

## References

Barrallo, J. (2011) Ovals and Ellipses in Architecture, *Proceedings of ISAMA 2011 Tenth Interdisciplinary Conference of the International Society of the Arts, Mathematics, and Architecture*, Columbia College, Chicago, Illinois, pp. 9-18.

Bazik, D., Dželebdžić, O. (2012) Web-based support of spatial planning in Serbia, *SPATIUM International Review*, No. 28, December 2012, pp.67-73

Boltyanski, V., Martini, H., Soltan, V. (1999), *Geometric methods and optimization problems*, Kluwer Acad. Publ., Dordrecht, Netherlands.

Brazil, M., Graham, R. L., Thomas, D. A., Zachariasen, M. (august 2013), *On the History of the Euclidean Steiner Tree Problem*, Archive for History of Exact Sciences, Springer Berlin Heidelberg.

*Crux Mathematicorum* - problem 2215, (1997), pp. 121-122. http://journals.cms.math.ca/CRUX/

Descartes, R. (1638), *Oeuvres de Descartes (publiées par Charles Adam et Paul Tannery 1897-1913)*, Tome II, Paris.

Drezner, Z, Hamacher, H. W. (Eds.) (2002). *Facility location: Applications and Theory*. Springer-Verlag, Berlin, pp. 1-11.

Erdös, P., Vincze, I. (1982) On the approximation of convex, closed plane curves by multifocal ellipses, *J. Appl. Probab.*, 19A, pp. 89–96

Farahani, R. Z., Hekmatfar, M. (2009) *Facility Location: Concepts, Models, Algorithms and Case Studies*, Chapter 1: Marzie Zarinbal, Distance Functions in Location Problems, Springer-Verlag Berlin Heidelberg, pp. 5-18 de Fermat P. (1679), *Oeuvres de Fermat*, Livre I, Paris, English edition 1891, Chapter V: Ad methodum de maxima et minima appendix , p. 153.

Groß, C., Strempel, T-K. (1998) On generalizations of conics and on a generalization of the Fermat-Torricelli problem, *American Mathematical Monthly*, pp. 732-743.

Harman P. M. (1990) *The Scientific Letters and Papers of James Clerk Maxwell*, Vol.1, 1846-1862, Cambridge University Press, pp. 35-62.

Heinrich, D. (1965) *100 Great Problems Of Elementary Mathematics their history and solution*, Dover Publication, Inc. New York, pp. 361-363.

Kirszenblat, D. (2011) *Dubins networks*, Thesis, Department of Mathematics and Statistics, The University of Melbourne.

Khilji, M. J. (2004) Multi Foci Closed Curves, *Journal of Theoretics* 6(6).

Kupitz, Y. S., Martini, H., Spirova, M. (2013) The Fermat–Torricelli Problem, Part I: A Discrete Gradient-Method Approach, *Journal of Optimization Theory and Applications*, pp. 1-23.

Mahon, B. (2003) *The Man Who Changed Everything – the Life of James Clerk Maxwell*, Hoboken, NJ: Wiley, p. 16.

Maxwell, J. C. (1847) On the description of oval curves, and those having a plurality of foci, *Proceedings of the Royal Society of Edinburgh*, Vol. ll, 1846., "On Oval Curves", 1847., "On Trifocal Curves", 1847.

Melzak, Z.A., Forsyth, J.S. (1977) Polyconics 1: Polyellipses and optimization, *Quart. Appl. Math.*, 35, pp. 239–255.

Mladenović, N. (2004) *Kontinualni lokacijski problemi*, Matematički institut SANU, Beograd.

Nie, J., Parrilo, P. A., Sturmfels, B. (2008) Semidefinite representation of the k-ellipse, *Algorithms in algebraic geometry*. Springer, New York, pp. 117-132.

Oikonomia, http://www.oikonomia.info/wp-content/uploads/2013/05/Southstream-route.jpg, accessed 24[th] Nov 2013.

de Pérgamo, A., Viviano, V. (1659) *De maximis et minimis geometrica divinatio in quintum conicorum Apollonio Pergaei*, Chapter: Ad lib de max et min, Apendix Monitvm, Probl.I, Prop VI, pp. 143-150.

Petrović, M. (2010) *Egg curves and generalisation Hugelschaffer's construction*, Magister thesis (in Serbian), Faculty of Architecture, University of Belgrade.

Petruševski, Lj. Dabić, M., Devetaković, M. (2010) Parametric curves and surfaces -MATHEMATICA demonstrations as a tool in exploration of architectural form, *SPATIUM International Review*, No. 22, July 2010, pp. 67-72.

Rosin, P. (2000). On Serlio's Construction of Ovals, *Nexus Network Journal*, vol. 2, no. 3.

Rosin, P. (2004) On the construction of ovals, *Proceedings of International Society for the Arts, Mathematics, and Architecture* (ISAMA), pp. 118-122.

Sahadevan, P.V. (1974) Evolution of 'N' polar curves by extension of focal points (by more than two) of an Ellipse, *Proceedings of the Indian Academy of Sciences* - Section A, 79(6), pp. 269-281.

Sahadevan, P.V. (1987) The theory of the egglipse - a new curve with three focal points, *Int. J. Math. Educ. Sci. Technol.* 18(1), pp. 29-39.

Sándor, J. (2005) On Certain Inequality by Visschers, *Octogon Math. Mag.* 13(2), pp. 1053-1054.

Sekino, J. (1999) n-Ellipses and the minimum distance sum problem, *The American mathematical monthly*, *106*(3), pp.193-202.

Simpson, T. (1750) *Doctrine and applications of fluxions*, London, p. 27.

Soothill, G. (2010) *The Euclidean Steiner Problem*, Report, Department of Mathematical Science, Durham University, England (Undergraduate Departmental Prizes).

Teodorović, D. (2009) *Transportne mreže*, Poglavlje 9: Lokacijski problem, Saobraćajni fakultet, Beograd, pp.389-399.

Torricelli E., Loria G., Vassura G. (1919) *Opere de Evangelista Torricelli* Vol I, Part 2, English edition, Faënza, pp. 90–97.

Varga, A., Vincze, Cs. (2008) On a lower and upper bound for the curvature of ellipses with more than two foci, *Expositiones Mathematicae* 26(1), pp.55–77.

Volek, J. (2006) Location analysis – Possibilities of use in public administration, *Verejná správa* 2006, Univerzita Pardubice, pp. 84–85.

Watanabe, D., Majima, T., Takadama, K., Katuhara, M. (2009) Generalized Weber Model for Hub Location of Air Cargo, The Eighth International Symposium on Operations Research and Its Applications (ISORA'09), Zhangjiajie, China, pp. 124–131.

Weber, A. (1929) [translated by Carl J. Friedrich from Weber's 1909 book *Über den Standort der Industrien*], *Alfred Weber's Theory of the Location of Industries*, Chicago, The University of Chicago Press.

Zuber, M., Štulić, R. (2006) O grafičkim mogućnostima programskog paketa Maple, The 23rd Conference on Descriptive Geometry and Engineering Graphics, MoNGeometrija 2006, Novi Sad Proceedings, pp. 50-59.